\documentclass [12pt] {amsart}
\usepackage{amsmath, amsfonts, latexsym, array, epsfig, amssymb, psfrag}
\usepackage[all]{xypic}
\usepackage{graphicx}
\newtheorem{thm}{Theorem}[section]
\newtheorem{cor}[thm]{Corollary}
\newtheorem{lem}[thm]{Lemma}
\newtheorem{prop}[thm]{Proposition}
\newtheorem{cla}[thm]{Claim}

\newtheorem{rem}[thm]{Remark}
\newtheorem{que}[thm]{Question}

\begin{document}
\title[Commuting Diffeomorphisms]
{Trivial centralizers for Axiom A diffeomorphisms}
\author{Todd Fisher}
\address{Department of Mathematics, Brigham Young University, Provo, UT 84602}
\email{tfisher@math.byu.edu}
\thanks{Supported in part by NSF Grant \#DMS0240049}

\date{September 1, 2006}
\commby{}

\begin{abstract}
We show there is a residual set of non-Anosov $C^{\infty}$ Axiom A diffeomorphisms with the no cycles property whose elements have trivial centralizer.  If $M$ is a surface and $2\leq r\leq \infty$, then we will show there exists an open and dense set of of $C^r$ Axiom A diffeomorphisms with the no cycles property whose elements have trivial centralizer.  Additionally, we examine commuting diffeomorphisms preserving a compact invariant set $\Lambda$ where $\Lambda$ is a hyperbolic chain recurrent class for one of the diffeomorphisms.  
\end{abstract}

\maketitle

\section{Introduction}

Inspired by Hilbert's address in 1900 Smale was asked for a list of problems for the $21^{\mathrm{st}}$ century.
Problem $12$ deals with the centralizer of a ``typical'' diffeomorphism.  For $f\in \mathrm{Diff}^r(M)$ (the set of $C^r$ diffeomorphisms from $M$ to $M$) the centralizer of $f$ is
$$Z(f)=\{g\in \mathrm{Diff}^r(M)\, |\, fg=gf\}.$$
Let $r\geq 1$, $M$ be a smooth, connected, compact, boundaryless manifold, and $$T=\{ f\in\mathrm{Diff}^r(M)\, |\, Z(f)\textrm{ is trivial}\}.$$  Smale asks the following question.

\begin{que}\label{q.1}
Is $T$ dense in $\mathrm{Diff}^r(M)$?
\end{que}

A number of people have worked on these and related problems including Kopell~\cite{Kop1}, Anderson~\cite{And1}, Palis and Yoccoz~\cite{PY1}~\cite{PY2}, Katok~\cite{K1}, Burslem~\cite{Bur}, Togawa~\cite{Tog}, and Bonatti, Crovisier, Vago, and Wilkinson~\cite{BCW}.

Palis and Yoccoz~\cite{PY1}  are able to answer Question~\ref{q.1} in the affirmative in the case of $C^{\infty}$ Axiom A diffeomorphisms with the added assumption of strong transversality.  In~\cite{PY1} Palis and Yoccoz ask if their results extend to the case of Axiom A diffeomorphisms with the no cycles property.

The first results of the present work extend the results of Palis and Yoccoz to the case of Axiom A diffeomorphisms with the no cycles property.  Let $\mathcal{A}^r(M)$ denote the set of $C^{r}$ Axiom A diffeomorphisms with the no cycles property that are non-Anosov. 
Let $\mathcal{A}^r_1(M)$ denote the subset of $\mathcal{A}^r(M)$ containing a periodic sink or source.

\begin{thm} \label{t.2} There is an open and dense subset of $\mathcal{A}^{\infty}_1(M)$ whose elements have a trivial centralizer.
\end{thm}

\begin{thm}\label{t.1} 
Let $\mathrm{dim}(M)\geq 3$.  Then there is a residual set of $\mathcal{A}^{\infty}(M)$ whose elements have trivial centralizer.
\end{thm}

We can reduce the requirement of $r=\infty$ for certain Axiom A diffeomorphisms of surfaces.  We note that in the following result we are able to include the Anosov diffeomorphisms

\begin{thm}\label{t.codimone}
If $M$ is a surface and $2\leq r\leq \infty$, then there exists an open and dense set of $C^r$ Axiom A diffeomorphisms with the no cycles property whose elements have trivial centralizer.
\end{thm}

A general technique in studying the centralizer is to examine properties of the centralizer on the ``basic pieces'' of the recurrent points.   In the hyperbolic case one often looks at what is called a hyperbolic chain recurrent class. 
 An {\it $\epsilon$-chain} from a point $x$ to a
point $y$ for a diffeomorphism $f$ is a sequence $\{x=x_0,...,x_n=y\}$ such
that the $d(f(x_{j-1}),x_j)<\epsilon$ for all $1\leq j\leq n$. The
\textit{chain recurrent set of f} is denoted $\mathcal{R}(f)$ and
defined by:
$$\mathcal{R}(f)=\{x\in M\, |\,\textrm{ there is an }\epsilon\textrm{-chain from }x
\textrm{ to }x\textrm{ for all }\epsilon>0\}.$$
For a point $x\in\mathcal{R}(f)$ the {\it chain recurrent class} of $x$ consists of all points $y\in\mathcal{R}(f)$ such that for all $\epsilon>0$ there is an $\epsilon$-chain from $x$ to $y$ and an $\epsilon$-chain from $y$ to $x$.

If Smale's question can be answered in the affirmative
one would hope that the following is true:  for a residual set of diffeomorphisms $\mathcal{G}$ if $\Lambda$ is a hyperbolic chain recurrent class for $f\in\mathcal{G}$ and $g\in Z(f)$ that there is a dichotomy, either $g|_{\Lambda}$ is the identity or $g|_{\Lambda}$ is a hyperbolic chain recurrent class.
In regards to this dichotomy we have the following result.

\begin{thm}\label{t.hypcommute} Suppose $f\in\mathrm{Diff}^r(M)$ for any $r\geq 1$, $g\in Z(f)$, $\Lambda$ is a mixing hyperbolic chain recurrent class of $f$ and a hyperbolic set for $g$.  Then $\Lambda$ is a locally maximal hyperbolic set for $g$.
\end{thm}

We briefly remark that
the above statement provides some context for an open problem
in the theory of higher rank symbolic actions: if an expansive
homeomorphism $f$
commutes with a transitive shift of finite type, must $f$ be
topologically
conjugate to a shift of finite type?~\cite{Nas}
Specifically, Theorem~\ref{t.hypcommute} is a smooth analog to this question.

\section{Background}

We now review some basic definitions and facts about hyperbolic sets and commuting diffeomorphisms.
We assume that all of our maps are diffeomorphisms of a manifold to itself.

A compact set $\Lambda$ invariant under the action of $f$ is {\it
hyperbolic} if there exists a splitting of the  tangent space
$T_{\Lambda}f=\mathbb{E}^u\oplus \mathbb{E}^s$ and positive
constants $C$ and $\lambda<1$ such that, for any point
$x\in\Lambda$ and any $n\in\mathbb{N}$,
$$
\begin{array}{llll}
\| Df_{x}^{n}v\|\leq C \lambda^{n}\| v\|,\textrm{ for }v\in
E^{s}_x \textrm{, and}\\
\| Df_{x}^{-n}v\|\leq C \lambda^{n}\| v\|,\textrm{ for }v\in
E^{u}_x. \end{array}
$$

For $\epsilon>0$ sufficiently small and $x\in \Lambda$ the
\textit{local stable and unstable manifolds} are respectively:
$$
\begin{array}{llll}
W_{\epsilon}^{s}(x,f)=\{ y\in M\, |\textrm{ for all }
n\in\mathbb{N}, d(f^{n}(x), f^{n}(y))\leq\epsilon\},\textrm{
and}\\
W_{\epsilon}^{u}(x,f)=\{ y\in M\, |\textrm{ for all }
n\in\mathbb{N}, d(f^{-n}(x), f^{-n}(y))\leq\epsilon\}.
\end{array}$$
The \textit{stable and unstable manifolds} are respectively:
$$
\begin{array}{llll}
W^s(x,f)=\bigcup_{n\geq 0}f^{-n}\left(
W_{\epsilon}^s(f^n(x),f)\right) ,\textrm{ and}\\
W^u(x,f)=\bigcup_{n\geq
0}f^{n}\left(W_{\epsilon}^u(f^{-n}(x),f)\right). \end{array}
$$
For a $C^r$ diffeomorphism the stable and unstable manifolds of a hyperbolic set are $C^r$
injectively immersed submanifolds.

A point $x$ is \textit{non-wandering} for a diffeomorphism $f$ if for any neighborhood $U$ of $x$ there exists an $n\in\mathbb{N}$ such that $f^n(U)\cap U\neq\emptyset$.  The set of non-wandering points is denoted $\mathrm{NW}(f)$.  A diffeomorphism $f$ is Axiom A if $\mathrm{NW}(f)$ is hyperbolic and is the closure of the periodic points.

A hyperbolic set is {\it locally maximal} if there exists a neighborhood $U$ of $\Lambda$ such that $\Lambda=\bigcap_{n\in\mathbb{Z}}f^n(U)$.
Locally maximal hyperbolic sets have some special properties.  First, we have the standard result called the
Shadowing Theorem, see~\cite[p. 415]{Rob}.  Let $\{x_j\}_{j=j_1}^{j_2}$
be an $\epsilon$-chain for $f$.  A point $y$ {\it
$\delta$-shadows} $\{x_j\}_{j=j_1}^{j_2}$ provided
$d(f^j(y),x_j)<\delta$ for $j_1\leq j\leq j_2$.

\begin{thm}(Shadowing Theorem) If $\Lambda$ is a locally maximal hyperbolic
set, then given any $\delta>0$ there exists an $\epsilon>0$ and $\eta>0$ such
that if $\{ x_j\}_{j=j_1}^{j_2}$ is an $\epsilon$-chain for $f$ with
$d(x_j, \Lambda)<\eta$, then there is a $y$ which $\delta$-shadows
$\{x_j\}_{j=j_1}^{j^2}$.  If the $\epsilon$-chain is periodic, then $y$ is
periodic.  If $j_2=-j_1=\infty$, then $y$ is unique and $y\in\Lambda$.
\end{thm}

The Shadowing Theorem also implies the following:

\begin{cor}\label{c.equivalent} If $\Lambda$ is a locally maximal hyperbolic set of a diffeomorphism $f$, then
$\mathrm{cl}(\mathrm{Per}(f|_{\Lambda}))=\mathrm{NW}(f|_{\Lambda})=\mathcal{R}(f|_{\Lambda})$.
\end{cor}

An additional consequence of the Shadowing Theorem is the structural stability of hyperbolic sets.  The following is a classical result, see~\cite[p.  571-572]{KH1}.

\begin{thm}(Structural stability of hyperbolic sets) Let $f\in \mathrm{Diff}(M)$ and $\Lambda$ be a hyperbolic set for $f$.  Then for any neighborhood $V$ of $\Lambda$ and every $\delta>0$ there exists a neighborhood $\mathcal{U}$ of $f$ in $\mathrm{Diff}(M)$ such that for any $g\in \mathcal{U}$ there is a hyperbolic set $\Lambda_g\subset V$ and a homeomorphism $h:\Lambda_g\rightarrow \Lambda$ with $d_{C^0}(\mathrm{id},h) + d_{C^0}(\mathrm{id},h^{-1})<\delta$ and $h\circ g|_{\Lambda_g}=f|_{\Lambda}\circ h$.  Moreover, $h$ is unique when $\delta$ is sufficiently small.
\end{thm}

If $X$ is a compact
set of a smooth manifold $M$ and $f$ is a continuous map from $M$
to itself, then $f|_X$ is {\it transitive} if for any open
sets $U$ and $V$ of $X$ there exists some $n\in\mathbb{N}$ such
that $f^n(U)\cap V\neq\emptyset$. A set $X$ is {\it mixing} if for
any open sets $U$ and $V$ in $X$ there exists an $N\in\mathbb{N}$
such that $f^n(U)\cap V\neq\emptyset$ for all $n\geq N$.

A standard result for locally maximal hyperbolic sets is the following Spectral Decomposition
Theorem~\cite[p. 575]{KH1}. 

\begin{thm}(Spectral Decomposition)\label{t.spectraldecomp}
Let $f\in\mathrm{Diff}^r(M)$  and $\Lambda$ a locally maximal hyperbolic set for $f$.  Then there
exist disjoint closed sets $\Lambda_1,...,\Lambda_m$ and a
permutation $\sigma$ of $\{1,...,m\}$ such that
$\mathrm{NW}(f|_{\Lambda})=\bigcup_{i=1}^m\Lambda_i$,
$f(\Lambda_i)=\Lambda_{\sigma(i)}$, and when $\sigma^k(i)=i$ then
$f^k|_{\Lambda_i}$ is topologically mixing.
\end{thm}

Corollary~\ref{c.equivalent} implies that Theorem~\ref{t.spectraldecomp} can be stated for a decomposition of the chain recurrent set of $f$ restricted to $\Lambda$ where $\Lambda$ is a locally maximal hyperbolic set. 

In the case where $f$ is Axiom A we have the following version of the Spectral Decomposition Theorem.
\\

\begin{thm}~\cite[p. 422]{Rob} Let $f\in\mathrm{Diff}^1(M)$ and assume that $f$ is Axiom A.  Then there are a finite number of sets $\Lambda_1,...,\Lambda_N$ closed, pairwise disjoint, and invariant by $f$ such that $\mathrm{NW}(f)=\bigcup_{i=1}^N\Lambda_i$.  Furthermore, each $\Lambda_i$ is topologically transitive.
\end{thm}

In the theorem above the sets $\Lambda_i$ are called \textit{basic sets}.  We define a relation $\ll$ on the basic sets $\Lambda_1,...,\Lambda_m$ given by the Spectral Decomposition Theorem as follows:  $\Lambda_i\ll\Lambda_j$ if $$(W^u(\Lambda_i)-\Lambda_i)\cap (W^s(\Lambda_j)-\Lambda_j)\neq \emptyset.$$  
A \textit{$k$-chain} is a sequence $\Lambda_{j_1},...,\Lambda_{j_k}$ where $\Lambda_{j_i}\neq \Lambda_{j_l}$ for $i,l\in[1,k]$ and $$\Lambda_{j_1}\ll\Lambda_{j_2}\ll ... \ll\Lambda_{j_k}.$$
A \textit{$k$-cycle} is a sequence of basic sets $\Lambda_{j_1},...,\Lambda_{j_k}$ such that 
$$\Lambda_{j_1}\ll\Lambda_{j_2}\ll ... \ll\Lambda_{j_k}\ll\Lambda_{j_1}.$$     A diffeomorphism is {\it Axiom A with the no cycles property} if the diffeomorphisms is Axiom A and there are no cycles between the basic sets given by the Spectral Decomposition Theorem.

The set $\mathcal{R}(f)$ is hyperbolic if and only if $f$ is Axiom A with the no cycles property.  For any $r\geq 1$ the set of $C^r$ diffeomorphisms with $\mathcal{R}(f)$ hyperbolic is open.

For a hyperbolic set $\Lambda$ let
$$W^s(\Lambda)=\{x\in M\,|\, \lim_{n\rightarrow\infty}d(f^n(x),\Lambda)=0\}.$$
If $\Lambda$ is a topologically transitive locally maximal hyperbolic set and $p\in\mathrm{Per}(f)\cap\Lambda$, then $W^s(\mathcal{O}(p))$ is dense in $W^s(\Lambda)$.  If $\Lambda$ is a mixing locally maximal hyperbolic set and $p\in\mathrm{Per}(f)\cap\Lambda$, then $W^s(p)$ is dense in $W^s(\Lambda)$.

The following proposition found in~\cite{ABCD} will be used in the proof of Theorem~\ref{t.hypcommute}.

\begin{prop}\label{p.abcd}
If $\Lambda$ is a hyperbolic chain recurrent class, then there exists a neighborhood $U$ of $\Lambda$ such that $\mathcal{R}(f)\cap U=\Lambda$.
\end{prop}

%Let $p$ be a hyperbolic periodic point for $f$.  A {\it homoclinic point} for $p$ is a point contained in $W^s(p)\cap W^u(p)$.  The {\it homoclinic class} of $p$ is the closure of all the homoclinic points.  A {\it hyperbolic homoclinic class} is a homoclinic class with a hyperbolic structure.

A set $X\subset M$ has an {\it attracting neighborhood} if there
exists a neighborhood $V$ of $X$ such that
$X=\bigcap_{n\in\mathbb{N}}f^n(V)$. A set $X\subset M$ has a {\it
repelling neighborhood} if there exists a neighborhood $U$ of $X$
such that $X=\bigcap_{n\in\mathbb{N}}f^{-n}(U)$. A set $\Lambda\subset M$
is called a {\it hyperbolic attractor} ({\it hyperbolic repeller})
if $\Lambda$ is a transitive hyperbolic set for a diffeomorphism
$f$ with an attracting neighborhood (a repelling neighborhood).  A
hyperbolic attractor (repeller) is {\it non-trivial} if it is not
the orbit of a periodic sink (source).

\begin{rem}\label{r.1}
For Axiom A diffeomorphisms with the no cycles property there is an open and dense set of points of the manifold that are in the basin of a hyperbolic attractor.
\end{rem}

We now review some basic properties of commuting diffeomorphisms.
Let $f$ and $g$ be commuting diffeomorphisms.  Let
$\mathrm{Per}^n(f)$ be the periodic points of period $n$ for $f$
and $\mathrm{Per}^n_h(f)$ denote the hyperbolic periodic points in
$\mathrm{Per}^n(f)$.  If $p\in\mathrm{Per}^n(f)$, then
$g(p)\in\mathrm{Per}^n(f)$ so $g$ permutes the points of
$\mathrm{Per}^n(f)$.  Furthermore, if
$p\in\mathrm{Per}^n(f)$, then
$$T_{g(p)}f^nT_pg=T_pgT_pf^n.$$
Hence, the linear maps $T_{g(p)}f^n$ and $T_pf^n$ are similar.
If $p\in\mathrm{Per}^n_h(f)$, then  $g(p)\in \mathrm{Per}^n_h(f)$.  Since $\#(\mathrm{Per}^n_h(f))<\infty$ it
follows that if $p\in\mathrm{Per}^n_h(f)$, then $g(p)\in\mathrm{Per}(g)$.  Additionally, If $p\in\mathrm{Per}^n_h(f)$, then
$$g(W^u(p,f))=W^u(g(p),f)\textrm{ and }g(W^s(p,f))=W^s(g(p),f).$$

\section{Trivial centralizer for Axiom A diffeomorphisms with no cycles}\label{s.axioma}

We will show that Theorems~\ref{t.2} and~\ref{t.1} will follow from extending the results in~\cite{PY1} if one can show the following theorem:

\begin{thm}\label{t.rigid} There exists and open and dense set $\mathcal{V}$ of $\mathcal{A}^{\infty}(M)$ such that if $f\in\mathcal{V}$ and $g_1,g_2\in Z(f)$ where $g_1=g_2$ on a non-empty open set of $M$, then $g_1=g_2$.
\end{thm}

Before proceeding with the proof of Theorem~\ref{t.rigid} we review a result of Anderson~\cite{And1}.  Let $f\in\mathrm{Diff}^{\infty}(\mathbb{R}^n)$ be a contraction.  Anderson shows that if $g_1,g_2\in Z(f)$ and $g_1=g_2$ on an open set of $\mathbb{R}^n$, then $g_1=g_2$ on all of $\mathbb{R}^n$.

In the proof of Theorem~\ref{t.rigid} it is sufficient to show that there exists an open set $\mathcal{V}$ of $\mathcal{A}^{\infty}(M)$ such that if $f\in\mathcal{V}$, $g\in Z(f)$, and there exists an open set $U$ of $M$ where $g|_U=\mathrm{id}|_U$, then $g=\mathrm{id}_M$. 

Now suppose that $f\in\mathcal{A}^{\infty}(M)$ and $g\in Z(f)$ where $g$ is the identity for a non-empty open set $U$ of $M$.  Then $U$ intersects the basin of either a hyperbolic attractor or repeller $\Lambda$ for $f$ in an open set denoted $U_{\Lambda}$.  Let $p\in\mathrm{Per}^n(\Lambda)$.  Then there exists a $i\in\mathbb{N}$ such that $W^s(f^i(p))\cap U\neq\emptyset$.
Since $g(W^s(f^i(p)))=W^s(q)$ for some periodic point $q$ and $g$ is the identity on $U$ we know that $g(f^i(p))=f^i(p)$.   The map $f^n$ restricted to $W^s(f^i(p))$ is a contraction that commutes with $g$ restricted to $W^s(f^i(p))$.  Hence, $g$ is the identity on $W^s(f^i(p))$ from Anderson's result.  The density of $W^s(\mathcal{O}(p))$ in $W^s(\Lambda)$ implies that $g$ is the identity on $W^s(\Lambda)$.  A similar argument holds for repellers.  

To prove Theorem~\ref{t.rigid} we then need a way to connect the basins of adjacent attractors so that if $g$ is the identity in one it will be the identity in the other.  To do this we will prove the next proposition.  We note that the next proposition is similar to the Lemma in the proof of Theorem 1 from~\cite[p. 85]{PY1}.

\begin{prop}\label{l.rigid} There exists an open and dense set  $\mathcal{V}$ of $\mathcal{A}^{r}(M)$, $1\leq r\leq \infty$, such that if $f\in\mathcal{V}$ and $\Lambda$ and $\Lambda'$ are attractors for $f$ where $$\overline{W^s(\Lambda)}\cap \overline{W^s(\Lambda')}\neq\emptyset,$$ then there exists a hyperbolic repeller $\Lambda_r$ such that $$W^s(\Lambda)\cap W^u(\Lambda_r)\neq\emptyset\textrm{ and }W^s(\Lambda')\cap W^u(\Lambda_r)\neq\emptyset.$$
\end{prop}

Before proving the above proposition we show how it implies Theorem~\ref{t.rigid}.
\\

\noindent\textbf{Proof of Theorem~\ref{t.rigid}.}  Let $\mathcal{V}$ be open and dense in $\mathcal{A}^{\infty}(M)$ satisfying Proposition~\ref{l.rigid} and let $f\in\mathcal{V}$.  Since $f\in\mathcal{A}^{\infty}(M)$ we know there is an open and dense set of $M$ contained in the basin of hyperbolic attractors.  Denote the hyperbolic attractors of $f$ as $\Lambda_1,...,\Lambda_k$.  Let $g\in Z(f)$ such that $g$ is the identity on a non-empty open set $U$ contained in $M$.  Then there exists some attractor $\Lambda_i$ where $1\leq i\leq k$ such that $W^s(\Lambda_i)\cap U\neq\emptyset$.   Hence, $g$ is the identity on $W^s(\Lambda_i)$.

For any attractor $\Lambda_j$ such that $$\overline{W^s(\Lambda_i)}\cap\overline{W^s(\Lambda_j)}\neq\emptyset$$ there exists a repeller $\Lambda_r$ such that  $$W^s(\Lambda_i)\cap W^u(\Lambda_r)\neq\emptyset\textrm{ and } W^s(\Lambda_j)\cap W^u(\Lambda_r)\neq\emptyset.$$  This follows from Proposition~\ref{l.rigid}.   It then follows that $g$ is the identity on $W^u(\Lambda_r)$ and $W^s(\Lambda_j)$ since the intersection of the basins for an attractor and a repeller is an open set.

Continuing the argument we see that $g$ is the identity on $\bigcup_{n=1}^k W^s(\Lambda_n)$.  Hence, $g$ is the identity on all of $M$ from Remark~\ref{r.1}.  $\Box$

We now state and prove two lemmas that will be helpful in proving Proposition~\ref{l.rigid}.  

\begin{lem}\label{l.transverse1} There exists an open and dense set $\mathcal{V}_1$ of $\mathcal{A}^{r}(M)$ for $1\leq r\leq \infty$ such that if $f\in\mathcal{V}_1$, $\Lambda$ is a hyperbolic repeller for $f$, and $\Lambda=\Lambda_0\ll \Lambda_1\ll ...\ll\Lambda_k$   , then $\Lambda\ll \Lambda_k$.
\end{lem}

\noindent{\bf Proof.}  Let $\mathcal{U}$ be a connected component of $\mathcal{A}^r(M)$.  Let $\Lambda_0,...,\Lambda_M$ be basic sets such that $\Lambda_0,...,\Lambda_j$ are the hyperbolic repellers for each $f\in\mathcal{U}$.

We will prove the lemma inductively on $k$.  For $k=1$ the statement is trivially true.  Assume for $k\geq 1$ that there is an open and dense set $\mathcal{U}_k$ of $\mathcal{U}$ such that if $\Lambda=\Lambda_{n_0}\ll...\ll\Lambda_{n_k}$ where $0\leq n_0\leq j$, then $\Lambda\ll \Lambda_{n_k}$.

Fix 
$$\mathbf{\alpha}=(\alpha_0,...,\alpha_{k+1})\in\{ 0,..., j\}\times\{ j+1,...,M\}^{k+1}.$$  Let $I$ be the set of all such $\alpha$ 
and let $f\in\mathcal{U}_k$ such that $\Lambda_{\alpha_0}\ll ... \ll\Lambda_{\alpha_{k+1}}$ for $f$.  Since $f\in\mathcal{U}_k$ we know that $\Lambda_{\alpha_0}\ll \Lambda_{\alpha_k}\ll\Lambda_{\alpha_{k+1}}$.

The next claim will show there is an arbitrarily small $C^r$ perturbation of $f$ such that $\Lambda_{\alpha_0}\ll \Lambda_{\alpha_{k+1}}$.

\begin{cla}\label{c.perturb}  If $f\in\mathcal{A}^r(M)$, $\Lambda\ll \Lambda_1\ll \Lambda_2$, $\Lambda$ is a repeller, $y\in W^u(\Lambda_1)\cap W^s(\Lambda_2)$, and $U$ is a sufficiently small neighborhood of $y$, then there exists an arbitrarily small $C^r$ perturbation $\tilde{f}$ of $f$ with support in $U$ such that $f(y)\in W^u(\Lambda)\cap W^s(\Lambda_2)$ for $\tilde{f}$.
\end{cla}

\noindent \textbf{Proof of Claim~\ref{c.perturb}.}  Let $p$ be a periodic point of $\Lambda_1$.  Since $W^s(\mathcal{O}(p))$ is dense in $W^s(\Lambda_1)$ and $W^u(\Lambda)$ is open we know that there exists some $m$ such that $W^u(\Lambda)\cap W^s(f^m(p))\neq\emptyset$.  Let $x\in W^u(\Lambda)\cap W^s(f^m(p))$.  If we take a transversal to $W^s(f^m(p))$ at $x$ such that the transversal is contained in $W^u(\Lambda)$, then the Inclination Lemma (or $\lambda$-lemma)~\cite[p. 122]{BS} implies that the transversal accumulates on $W^u(f^m(p))$.   By the invariance of $W^u(\Lambda)$ the same holds for any power of $p$.

Let $y\in W^u(\Lambda_1)\cap W^s(\Lambda_2)$.  Then there exists an $n$ such that $y\in \overline{W^u(f^n(p))}$ and hence $W^u(\Lambda)$ accumulates on $y$. 

Since $y$ is a wandering point there exists a sufficiently small neighborhood $U$ of $y$ such that $f^n(U)\cap U=\emptyset$  for all $n\in\mathbb{Z}-\{0\}$, and  $U$ is disjoint from a neighborhood of $\Lambda\cup \Lambda_1\cup\Lambda_2$.

Let $y_k$ be a sequence in $W^u(\Lambda)$ converging to $y$.  Then there exists an arbitrarily small $C^r$ perturbation $\tilde{f}$, with support in $U$, of $f$ such that $y_k$ gets mapped to $f(y)$ for some $k$ sufficiently large.  We know that $\tilde{f}^{-n}(y_k)=f^{-n}(y_k)$  and $\tilde{f}^n(y)=f^n(y)$ for all for all $n\in\mathbb{N}$.  Hence $f(y)\in W^u(\Lambda)\cap W^s(\Lambda_2)$ for $\tilde{f}$.
$\Box$
\\

We now return to the proof of the lemma.
The previous claim shows that by an arbitrarily small perturbation $\tilde{f}$ we have 
$\Lambda_{\alpha_0}\ll \Lambda_{\alpha_{k+1}}$.  Since $\mathcal{U}_k$ is open we may assume $\tilde{f}\in\mathcal{U}_k$.  Since $W^u(\Lambda_{\alpha_0})$ is open and varies continuously with $\tilde{f}$ as does $W^s(\Lambda_{\alpha_{k+1}})$ we know that it is an open condition that $\Lambda_{\alpha_0}\ll ... \ll\Lambda_{\alpha_{k+1}}$ and $\Lambda_{\alpha_0}\ll \Lambda_{\alpha_{k+1}}$.  

Let 
$\mathcal{U}^0_{k,\mathbf{\alpha}}\subset \mathcal{U}_k$ such that for all $f\in\mathcal{U}^0_{k,\mathbf{\alpha}}$ there exists some $\epsilon>0$ where $\Lambda_{\alpha_0}\ll ... \ll\Lambda_{\alpha_{k+1}}$ is not a chain for all $g\in B_{\epsilon}(f)$.  Let $\mathcal{U}^1_{k,\mathbf{\alpha}}\subset \mathcal{U}_k$ such that $\Lambda_{\alpha_0}\ll ... \ll\Lambda_{\alpha_{k+1}}$ and $\Lambda_{\alpha_0}\ll \Lambda_{\alpha_{k+1}}$.  From the previous argument we know that $\mathcal{U}^1_{k,\mathbf{\alpha}}$ is open.  We want to show that $\mathcal{U}^1_{k,\mathbf{\alpha}}$ is dense in $\mathcal{U}_k-\mathcal{U}^0_{k,\mathbf{\alpha}}$.

Let $f\in \mathcal{U}_k-\mathcal{U}^0_{k,\mathbf{\alpha}}$.  Then there exists a sequence of $f_n$ converging to $f$ such that for each $f_n$ we have $\Lambda_{\alpha_0}\ll ...\ll\Lambda_{\alpha_{k+1}}$.  Hence, there exists a sequence $\tilde{f}_n$ converging to $f$ such that each $\tilde{f_n}\in \mathcal{U}^1_{k,\mathbf{\alpha}}$.

Define $\mathcal{U}_{k,\alpha}=\mathcal{U}^0_{k,\mathbf{\alpha}}\cup\mathcal{U}^1_{k,\mathbf{\alpha}}$.  The set $\mathcal{U}_{k,\alpha}$ is then open and dense in $\mathcal{U}_k$.  Define
$$\mathcal{U}_{k+1}=\bigcap_{\alpha\in I}\mathcal{U}_{k,\alpha}.$$  Since the set $I$ is finite we know that $\mathcal{U}_{k+1}$ is open and dense in $\mathcal{U}_k$.
The no cycles property implies there are no chains of length $M+1$.  Define $\mathcal{V}_1=\mathcal{U}_M$ this will be open and dense in $\mathcal{A}^r(M)$ and by construction if $\Lambda$ is a hyperbolic repeller for $f\in\mathcal{V}_1$, and $\Lambda=\Lambda_0\ll \Lambda_1\ll ...\ll\Lambda_k$   , then $\Lambda\ll \Lambda_k$ for $f$.
  $\Box$

\begin{lem}\label{l.transverse2} There exists an open and dense set $\mathcal{V}_2$ of $\mathcal{A}^r(M)$ for $1\leq r\leq\infty$ such that if $f\in\mathcal{V}_2$, $\Lambda$ is a hyperbolic attractor for $f$, $\Lambda'$ is a basic set for $f$, $\Lambda_r$ is a hyperbolic repeller for $f$ with $\Lambda_r\ll \Lambda'$, and $\overline{W^s(\Lambda)}\cap W^u(\Lambda')\neq \emptyset$, then $\Lambda_r\ll \Lambda$.
\end{lem}

\noindent{\bf Proof.}  Let $\mathcal{C}$ be a connected component of $\mathcal{A}^r(M)\cap\mathcal{V}_1$ where $\mathcal{V}_1$ is an open and dense set of diffeomorphisms satisfying Lemma~\ref{l.transverse1} and let $\Lambda_0,...,\Lambda_M$ be the basic sets for $f\in\mathcal{C}$ where  $\Lambda_0,...,\Lambda_j$ are the hyperbolic attractors and $\Lambda_J,\Lambda_{J+1}...,\Lambda_M$ are hyperbolic repellers.   Before proceeding with the proof of the lemma we prove a claim.

\begin{cla}\label{c.attractrepell}
If $f\in\mathcal{C}$ satisfies the following:
\begin{itemize}
\item $\Lambda_a$ is a hyperbolic attractor for $f$, 
\item $\Lambda_b$ is a basic set for $f$, 
\item $\Lambda_r$ is a repeller for $f$, 
\item $\overline{W^s(\Lambda)}\cap W^u(\Lambda')\neq\emptyset$, and 
\item $\Lambda_r\ll \Lambda_b$,
\end{itemize} 
then there exists an arbitrarily small perturbation of $f$ such that $\Lambda_r\ll \Lambda_a$. 
\end{cla}

\noindent{\bf Proof of claim.}
Lemmas 2.4 and 2.5 in~\cite{Shub1} show that 
$$\overline{W^s(\Lambda_a)}\cap \Lambda_b\neq \emptyset\textrm{ and }\overline{W^s(\Lambda_a)}\cap (W^u(\Lambda_b)-\Lambda_b)\neq \emptyset.$$  Let $x\in \overline{W^s(\Lambda_a)}\cap (W^u(\Lambda_b)-\Lambda_b)$.   Since $x$ is wandering there exists a neighborhood $V$ of $x$ such that $f^n(V)\cap V=\emptyset$ for all $n\in\mathbb{Z}-\{0\}$.

Since $\Lambda_r\ll \Lambda_b$ we know that $W^u(x)\subset \overline{W^u(\Lambda_r)}$.
As in the proof of Claim~\ref{c.perturb} there exists a $C^r$ small perturbation $\tilde{f}$ with support in $V$ such that $f(x)\in W^u(\Lambda_r)$ for $\tilde{f}$.

Since $f^n(V)\cap V=\emptyset$ for all $n\in\mathbb{N}$, the perturbation had support in $V$, and $f(x)\in W^u(\Lambda_r)\cap \overline{W^s(\Lambda_a)}$ for $\tilde{f}$ we know that $\Lambda_r\ll \Lambda_a$ for $\tilde{f}$. $\Box$

We now return to the proof of the lemma.
Let $$\alpha=(\alpha_1,\alpha_2, \alpha_3)\in\{0,...,j\}\times \{j+1,...,M\}\times\{J,..., M\}$$ and $I$ be the set of all such $\alpha$.  %Let $\mathcal{C}^0_{\alpha}$ be the set of $f\in\mathcal{C}$ such that there is a neighborhood $\mathcal{V}$ of $f$ in $\mathcal{C}$ such that for each $g\in\mathcal{V}$ the sets $\overline{W^s(\Lambda_{\alpha_1})}$ and $W^u(\Lambda_{\alpha_2})$ do not intersect.   
Let $\mathcal{C}^0_{\alpha}$ be the set of all $f\in\mathcal{C}$ such that if $\overline{W^s(\Lambda_{\alpha_1})}\cap W^u(\Lambda_{\alpha_2})\neq\emptyset$ and $\Lambda_{\alpha_3}\ll \Lambda_{\alpha_2}$, then $\Lambda_{\alpha_3}\ll \Lambda_{\alpha_1}$.  We want to show that $\mathcal{C}_{\alpha}=\mathrm{int}\,\mathcal{C}^0_{\alpha}$ is dense in $\mathcal{C}$.  

Let 
$f\in\mathcal{C}-\mathcal{C}_{\alpha}$.    Then for all neighborhoods $U$ of $f$ there exists a function $g\in U$ satisfying the following:
\begin{itemize}
\item  $\overline{W^s(\Lambda_{\alpha_1})}\cap W^u(\Lambda_{\alpha_2})\neq\emptyset$, 
\item $\Lambda_{\alpha_3}\ll\Lambda_{\alpha_2}$, and 
\item $W^u(\Lambda_{\alpha_3})\cap W^s(\Lambda_{\alpha_1})=\emptyset$.  
\end{itemize}
Then from Claim ~\ref{c.attractrepell} there exists an arbitrarily small perturbation $\tilde{g}$ of $g$ such that $\Lambda_{\alpha_3}\ll \Lambda_{\alpha_1}$ for $\tilde{g}$.  Since the intersection of the basins for attractors and repellers is an open condition among the diffeomorphisms we know that $\tilde{g}\in \mathcal{C}_{\alpha}$.  Hence, $f\in\overline{\mathcal{C}_{\alpha}}$ and $\mathcal{C}_{\alpha}$ is open and dense in $\mathcal{C}$.

Let 
$$\mathcal{V}_2=\bigcap_{\alpha\in I} \mathcal{C}_{\alpha}.$$
Since $I$ is a finite set this will be an open and dense set in $\mathcal{C}$.
 $\Box$
 \\
 
 \noindent{\bf Proof of Proposition~\ref{l.rigid}.}  Let $\mathcal{V}=\mathcal{V}_1\cap \mathcal{V}_2$ where $\mathcal{V}_1$ and $\mathcal{V}_2$ are open and dense sets in $\mathcal{A}^r(M)$ satisfying Lemma~\ref{l.transverse1} and Lemma~\ref{l.transverse2}, respectively.  Let $f\in\mathcal{V}$ and $\Lambda$ and $\Lambda'$ be attractors such that $\overline{W^s(\Lambda)}\cap \overline{W^s(\Lambda')}\neq\emptyset$.  Fix $x\in \overline{W^s(\Lambda)}\cap \overline{W^s(\Lambda')}$.  Then $x\in W^u(\tilde{\Lambda})$ for some basic set $\tilde{\Lambda}$.  Then from Lemma~\ref{l.transverse2} there exists a hyperbolic repeller $\Lambda_r$ such that $\Lambda_r\ll \Lambda$ and $\Lambda_r\ll\Lambda'$.  $\Box$
\\

To extend the proofs of Theorems 2 and 3 from~\cite{PY1} to Theorems~\ref{t.2} and~\ref{t.1} we now need to show that the lack of strong transversality is not essential in the arguments.  

Let $\mathcal{U}(M)$ be the set of $C^{\infty}$ Axiom A diffeomorphisms with the strong transversality condition.  Let $\mathcal{U}_1(M)$ consist of all elements of $\mathcal{U}(M)$ that have a periodic sink or source.

To prove Theorem 2 in~\cite{PY1} it is shown there is an open and dense set $\mathcal{C}_1(M)$ of $\mathcal{U}_1(M)$ such that if $f\in\mathcal{C}_1(M)$, then there is a periodic sink (or source) $p$ such that if  $g\in Z(f)$, then $g=f^k$ in $W^s(p)$ ($W^u(p)$).  Theorem 1 in~\cite{PY1} (that is similar to Theorem~\ref{t.rigid} in the present work) is then used to connect the regions to show that $g$ is a power of $f$ for all of $M$.

Similarly, to prove Theorem 3 in~\cite{PY1} it is shown there is a set $\mathcal{C}(M)$ that is residual in $\mathcal{U}(M)$ if $\mathrm{dim}(M)\geq 3$ such that for any $f\in\mathcal{C}(M)$ there is a hyperbolic attractor (or repeller) $\Lambda$ for $f$ such that if  $g\in Z(f)$, then $g=f^k$ in $W^s(\Lambda)$ ($W^u(\Lambda)$).  Theorem 1 in~\cite{PY1} is then used to connect the regions to show that $g$ is a power of $f$ for all of $M$.

Let $p$ be a periodic point of period $k$ for a diffeomorphism $f\in\mathcal{A}^{\infty}(M)$ of an $n$-dimensional manifold. The periodic point $p$ is {\it non-resonant} if the eigenvalues of $Df^kp$ are distinct and for all $(j_1,...,j_n)\in\mathbb{N}^n$ such that $\sum j_k\geq 2$, we have
$$\lambda_i\neq \lambda_1^{j_1}...\lambda_n^{j_n}\textrm{ for all }1\leq i\leq n.$$
For a hyperbolic periodic point it is clear that non-resonance is an open condition.

Fix $p$ a non-resonant hyperbolic periodic point for $f$.  Let $n_1$ and $n_2$ be the dimensions of the stable and unstable manifolds of $p$, respectively.  Then there exists~\cite[p.\,90-91]{PY1} immersions $\mathcal{H}^s_p(f)$ and $\mathcal{H}^u_p(f)$ such that:
\begin{enumerate}
\item $\mathcal{H}^s_p(f)(\mathbb{R}^{n_1})=W^s(p)$ and $\mathcal{H}^u_p(\mathbb{R}^{n_2})(f)=W^u(p)$, 
\item the immersions vary continuously with $f$,
\item $A^s_p(f)=\mathcal{H}^s_p(f)^{-1}\circ f^k|_{W^s(p)}\circ \mathcal{H}^s_p(f)$ is a non-resonant linear contraction of $\mathbb{R}^{n_1}$,  and
\item $A^u_p(f)=\mathcal{H}^u_p(f)^{-1}\circ f^{-k}|_{W^u(p)}\circ \mathcal{H}^u_p(f)$ is a non-resonant linear contraction of $\mathbb{R}^{n_2}$.
\end{enumerate}

Let $\mathcal{J}_p(f)=(W^s(p)\cap W^u(p))-\{p\}$.  Define a map $\varphi_p$ from $\mathcal{J}_p(f)$ into $\mathbb{R}^n$ by
$$\varphi_p(q)=(\mathcal{H}^s_p(f)^{-1}(q), \mathcal{H}^u_p(f)^{-1}(q)).$$

Let $\bar{\mathcal{J}}_p(f)=\varphi(\mathcal{J}_p(f))$.  For $f\in\mathcal{U}(M)$ the transversality of the stable and unstable manifolds implies that $\bar{\mathcal{J}}_p(f)$ is discrete and closed in $\mathbb{R}^n$.   In the proofs of Theorems 2 and 3 in~\cite{PY1} this is only one place where the transversality is essential.  In the case of Axiom A diffeomorphisms with the no cycles property it will not be true,  generally, that $\bar{\mathcal{J}}_p(f)$ is discrete or closed.

However, in the proofs of Theorems~\ref{t.2} and~\ref{t.1} the set $\bar{\mathcal{J}}_p(f)$ is only needed for $p$ a periodic point of a hyperbolic attractor.  In this case we know that $W^u(p)$ is contained in the attractor and the hyperbolic splitting says that $$(W^s(p)\cap W^u(p))-\{p\}=(W^s(p)\pitchfork W^u(p))-\{p\}.$$  Hence, the arguments for Theorems 2 and 3 in~\cite{PY1} carry over to the case of Axiom A diffeomorphisms with the no cycles property.  This then shows the following two theorems.

\begin{thm}
There is an open and dense set $\mathcal{C}_1(M)$ of $\mathcal{A}^{\infty}_1(M)$ such that if $f\in\mathcal{C}_1(M)$ there is a periodic sink or source $p$ such that if $g\in Z(f)$, then $g=f^k$ in $W^s(p)$ or $W^u(p)$. 
\end{thm}

\begin{thm}
There is a set $\mathcal{C}(M)$ that is residual in $\mathcal{A}^{\infty}(M)$ if $\mathrm{dim}(M)\geq 3$ such that for any $f\in\mathcal{V}(M)$ there is a hyperbolic attractor or repeller $\Lambda$ for $f$ such that if $g\in Z(f)$, then $g=f^k$ in $W^s(\Lambda)$ or $W^u(\Lambda)$.
\end{thm}

The proofs of Theorems~\ref{t.2} and~\ref{t.1} now follow from Theorem~\ref{t.rigid} and the above two theorems.

\section{Axiom A diffeomorphisms of surfaces}

In this section we prove Theorem~\ref{t.codimone}.  This Theorem can be seen as an extension of Theorem 2 in~\cite{Roc} where Rocha examines $C^1$ centralizers of $C^{\infty}$ Axiom A diffeomorphisms of surfaces.  The important difference is Proposition~\ref{l.rigid}, which allows connections for the basins of attractors and repellers in the $C^r$ setting where $2\leq r\leq \infty$.

Throughout this section assume that $M$ is a compact surface.  Let $\mathcal{A}_0^r(M)$ be the set of all Axiom A diffeomorphisms of $M$ with the no cycles property.  Notice that $\mathcal{A}_0^r(M)$ contains all Anosov diffeomorphisms.  Let $\mathcal{V}_0$ be an open and dense set of $\mathcal{A}_0^r(M)$, where $2\leq r\leq\infty$, satisfying Proposition~\ref{l.rigid}.
Let $\mathcal{C}$ be a connected component of $\mathcal{V}_0$ and fix $N\in\mathbb{N}$ such that for each $f\in\mathcal{C}$ and $\Lambda$ a basic set for $f$ there is a periodic point in $\Lambda$ of period $k$ where $k\leq N$.

Let $\mathcal{C}_N$ be an open and dense set of $\mathcal{C}$ satisfying the following:
\begin{enumerate}
\item  If $f\in\mathcal{C}_N$ and $p$ and $p'$ are periodic points of period $k\leq N$ for $f$, then $T_pf^k$ and $T_{p'}f^k$ are conjugate if and only if $p'$ is in the orbit of $p$.
\item If $f\in \mathcal{C}_N$ and $p$ is a periodic point of period $k\leq N$ for $f$, then the eigenvalues of $p$ are non-resonant.
\end{enumerate}

Before proceeding with the proof of Theorem~\ref{t.codimone} we show there is an open and dense set $\mathcal{V}_1$ of $\mathcal{C}_N$ such that for each $f\in\mathcal{V}_1$ if $p$ a periodic point of period $k\leq N$ and $g\in Z(f)$, then there is a linearization of $W^s(p)$ and $W^u(p)$ for $f$ that will also be a linearization for $g$.

We first show this for saddle periodic points.
We will use the following theorem of  Sternberg from~\cite{Ste}.

\begin{thm}\label{t.c1reals}
If $g\in\mathrm{Diff}^2(\mathbb{R})$, $g(0)=0$, and $g'(0)=a$ where $|a| \neq 1$, then there exists a $C^1$ map $\varphi:\mathbb{R}\rightarrow\mathbb{R}$ with a $C^1$ inverse such that $\varphi \circ g\circ \varphi^{-1}(x)=ax$ for all $x$ sufficiently small.
\end{thm}

In the case where the contraction has basin of attraction the entire real line we have the following corollary that will be useful in proving Theorem~\ref{t.codimone}.

\begin{cor}\label{c.c1reals}
If $g\in\mathrm{Diff}^2(\mathbb{R})$ is a contraction of the reals fixing the origin with $g'(0)=a$ where $|a|\neq 1$, then there exists a $C^1$ map $\varphi$ with $C^1$ inverse such that $\varphi\circ g\circ \varphi^{-1}(x)=ax$ for all $x$.
\end{cor}

Let $f\in\mathcal{C}_N$ and $p$ be a periodic saddle point of period $k\leq N$.  From the Stable Manifold Theorem there exist $C^2$ immersions $\psi_s:\mathbb{R}\rightarrow W^s(p)$ and $\psi_u:\mathbb{R}\rightarrow W^u(p)$.  Furthermore, the immersions vary continuously with $f$ in $\mathcal{C}_N$ and the maps $$\Psi_s=\psi^{-1}\circ f^k\circ\psi\textrm{ and }\Psi_u=\psi^{-1}\circ f^{-k}\circ \psi$$ are $C^2$ contractions of the reals with $\Psi'_s(0)$ and $\Psi'_u(0)$ both less than one in absolute value.

From Corollary~\ref{c.c1reals} there exist maps $\varphi_s$ and $\varphi_u$ such that 
$$F_s^p=\varphi_s\circ \Psi_s\circ \varphi_s^{-1}\textrm{ and }F_u^p=\varphi_u\circ \Psi_u\circ \varphi_u^{-1}$$ are $C^1$ linear contractions of the reals.

Let $g\in Z(f)$.  Since $g(p)\in\mathcal{O}(p)$ there is a unique $0\leq i<k$ such that $(g\circ f^i)(p)=p$.  Let $\bar{g}=g\circ f^i$.  Define 
$$\begin{array}{llll}
G^p_s=\varphi_s\circ \psi^{-1}\circ \bar{g}^k\circ\psi  \circ \varphi_s^{-1}, \textrm{ and}\\
G^p_u=\varphi_u\circ \psi^{-1}\circ \bar{g}^k\circ\psi  \circ \varphi_u^{-1}.
\end{array}
$$
By definition $G^p_s,G^p_u\in\mathrm{Diff}^1(\mathbb{R})$.  Furthermore,
$$
\begin{array}{llll}
G^p_sF^p_s & = \varphi_s\circ \psi^{-1}\circ \bar{g}^kf\circ\psi  \circ \varphi_s^{-1} \\
 & =\varphi_u\circ \psi^{-1}\circ f\bar{g}^k\circ\psi  \circ \varphi_u^{-1}=  F^p_sG^p_s.
\end{array}
$$
The next lemma shows that $G^p_s$ and $G^p_u$ are linear maps of the reals fixing the origin.

\begin{lem}\label{l.commutereals}~\cite[p. 61]{KH1}  Any $C^1$ map defined on a neighborhood of the origin on the real line and commuting with a linear contraction $L:x\rightarrow \lambda x$, $|\lambda |<1$, is linear.
\end{lem}

For $p$ a periodic sink (or source) we use a different, but similar approach.  Let $\lambda_1$ and $\lambda_2$ be the eigenvalues for $p$, these are non-resonant, and $|\lambda_1|<|\lambda_2|<1$.  First, Hartman shows in~\cite{Hart} that there is a $C^1$ linearization of the basin of attraction for the periodic sink (or source) $p$.  Then Rocha shows in~\cite{Roc} that if $g\in Z(f)$ we don't know that $g$ is linearized, but $G_s^p=(\mu_1x,\mu_2y) + (0,h(x))$ where $h:\mathbb{R}\rightarrow\mathbb{R}$ is $C^1$ and satisfies
\begin{enumerate}
\item $h(0)=0$ and $h'(0)=0$, and
\item $h(\lambda_1x)=\lambda_2h(x)$ for all $x\in\mathbb{R}$.
\end{enumerate}

Let $\mathcal{V}_1\subset \mathcal{C}_N$ satisfying
\begin{enumerate}
\item[(A1)] if $p$ is a sink (source), then
\begin{enumerate}
\item each connected component of $W^{ss}(p)$ ($W^{uu}(p)$) is contained in the basin of some repeller, or
\item there is a saddle-like basic set $\Lambda_1$ and $x\in\Lambda_1$ such that $W^{ss}(p)$ and $W^u(x)$ ($W^{uu}(p)$ and $W^s(x)$) have a common point of transversal intersection, and
\end{enumerate}
\item[(A2)]
\begin{enumerate}
\item if $p$ is a sink or saddle and $q$ is a source such that $W^{ss}(p)\cap W^u(q)\neq\emptyset$, then $W^{ss}(p)$ and the strong unstable foliation of $W^u(q)$, denoted $\mathcal{F}^{uu}(q)$, have a common point of transversal intersection,
\item if $p$ is a source or a saddle and $q$ is a sink such that $W^{uu}(p)\cap W^s(q)\neq\emptyset$, then $W^{uu}(p)$ and the strong unstable foliation of $W^s(q)$, denoted $\mathcal{F}^{ss}(q)$, have a common point of transversal intersection.
\end{enumerate}
\end{enumerate}

It is clear that $\mathcal{V}_1$ is open and dense in $\mathcal{C}_N$.

\begin{thm}\label{t.linking}
There is an open and dense set $\mathcal{V}_2$ of $\mathcal{A}_0^r(M)$, where $2\leq r\leq\infty$, such that if $f\in\mathcal{V}_2$ and $g_1,g_2\in Z(f)$ satisfy  $g_1|_U=g_2|_U$ for some open set $U$ of $M$, then $g_1=g_2$.
\end{thm}

\noindent{\bf Proof.}  It is sufficient to look at an open and dense set $\mathcal{V}_1$ of $\mathcal{C}_N$ satisfying (A1) and (A2) above.  Let $f\in \mathcal{V}_1$ and suppose that $g\in Z(f)$ such that $g|_U=\mathrm{id}|_U$.  Suppose $U$ is contained in the basin of an attractor $\Lambda$.  

If $\Lambda$ is nontrivial, then the linearization of the stable manifold and the density of the stable manifold in the basin of attraction implies that $g|_{W^s(\Lambda)}=\mathrm{id}|_{W^s(\Lambda)}$.
If $\Lambda$ is the orbit of a sink $p$, then Rocha shows in~\cite{Roc} that $g|_{W^s(\mathcal{O}(p))}=
\mathrm{id}|_{W^s(\mathcal{O}(p))}$.  The theorem now follows since $\mathcal{V}_1$ is a subset of $\mathcal{V}_0$. $\Box$
\\

\noindent{\bf Proof of Theorem~\ref{t.codimone}}.
We first assume there is an attractor of saddle type.  Then the argument is an extension of one used in Theorem 3 part (a) of~\cite{PY1}.
We will give an explanation of the details. (We will follow the notation used in~\cite{PY1}.) Let $\Lambda$ be an attractor of saddle type.

Let $Z=\mathbb{R}^2\times (\mathbb{Z}/2\mathbb{Z})^2$ and $\Sigma$ be the hyperplane in $\mathbb{R}^2$ such that $\theta_1 +\theta_2=0$.  Define $\xi: Z\rightarrow\Sigma$ by 
$$\chi(\theta_1,\theta_2,\epsilon_1,\epsilon_2)=(\theta'_1,\theta'_2)$$
where $\theta'_i=\theta_i - \frac{1}{2}(\theta_1 + \theta_2)$.  Let 
$$\mathcal{L}=\{(\theta_1,\theta_2,\epsilon_1,\epsilon_2)\in Z\, |\, \theta_1=\theta_2=1\}$$
and for $\epsilon\in\mathcal{L}$ denote the cyclic subgroup generated by $\epsilon$ as $(\epsilon)$.

Let $Z_0=Z/(\epsilon)$ and $Z_1=\mathrm{ker}(\chi)/(\epsilon)$.

Let $D$ be the set of diagonal matrices in $\mathrm{GL}(2,\mathbb{R})$.  For $B\in D$ with diagonal entries $\lambda_1$ and $\lambda_2$ there is an isomorphism of $D$ onto $Z$ given by
$$\Theta_B(\mathrm{diag}(\mu_1,\mu_2))=(\frac{\mathrm{Log} |\mu_1|}{\mathrm{Log} |\lambda_1|},
 \frac{\mathrm{Log} |\mu_2|}{\mathrm{Log} |\lambda_2|}, \frac{\mu_1}{|\mu_1|},   \frac{\mu_2}{|\mu_2|} ).$$
 
 Let
 $$
 F^p=
 \begin{bmatrix}
F^p_s & 0\\
0 & (F^p_u)^{-1}
\end{bmatrix} \textrm{ and }
 G^p=
 \begin{bmatrix}
G^p_s & 0\\
0 & (G^p_u)^{-1}
\end{bmatrix}
$$

Then the isomorphism $\Theta_{F^p}$ sends $F^p$ to $$\epsilon=(1,1,\mathrm{sgn}(F^p_s),\mathrm{sgn}((F^p_u)^{-1}).$$  This then defines an isomorphism from $D/F^p$ to $Z_0$.  Denote $\bar{G}^p$ as the projection of $G^p$ in $Z_0$.
Define $\bar{\mathcal{J}}_p(f)$ as in Section~\ref{s.axioma}.

If $\Lambda$ is not all of $M$ (the non-Anosov case), then the next proposition, which is a restatement of Proposition 1 part (a) and (b) in~\cite{PY1} will show that the centralizer is trivial.

\begin{prop}
There is an open and dense set $\mathcal{V}$ of $\mathcal{C}_N$ such that if $f\in\mathcal{V}$ and $g\in Z(f)$, then 
\begin{enumerate}
\item no $\bar{G}\in Z_1$, $\bar{G}\neq 1_{Z_1}$ leaves $\bar{\mathcal{J}}_p(f)$ invariant, and
\item no $\bar{G}\in Z_0-Z_1$ leaves $\bar{\mathcal{J}}_p(f)$ invariant.
\end{enumerate}
\end{prop}

For the proof see the proof of Proposition 1 part (a) and (b) in~\cite{PY1}.

Then for any $f\in\mathcal{V}$ if $\Lambda$ is a non-trivial hyperbolic attractor or repeller and $g\in Z(f)$, then $g=f^k$ for some $k\in\mathbb{N}$ in $W^s(\Lambda)$ or $W^u(\Lambda)$ from the above proposition.  Theorem~\ref{t.linking} then shows that $g$ is a power of $f$ on all of $M$.

If $\Lambda=M$, then the result follows from Proposition 1 of~\cite{PY2}.
The only case left is where all the attractors and repellers are trivial.  This follows from the proof of Theorem 2 in~\cite{Roc}.   $\Box$

\section{Hyperbolic sets for commuting diffeomorphisms}

Before proceeding to the proof of Theorem~\ref{t.hypcommute} we review some results
from~\cite{Fis1}. Let $\Lambda$ be a hyperbolic set and $V$ be a neighborhood of $\Lambda$.  Theorem 1.5 in~\cite{Fis1} shows there is a
hyperbolic set $\tilde{\Lambda}\supset\Lambda$ contained in
$\Lambda_V=\bigcap_{n\in\mathbb{Z}}f^n(\overline{V})$ with a Markov partition.

Let $\Lambda$ be a hyperbolic set. In the proof of Theorem 1.5
in~\cite{Fis1} it is shown there exists a constant $\nu>0$ such
that for any $\nu$-dense set $\{p_i\}_{i=1}^N$ in $\Lambda$ there
is a subshift of finite type $\Sigma_A$ associated with the
transition matrix $A$ with
$$a_{ij}=\left\{
\begin{array}{llll}1 & \textrm{if} & d(f(p_i),p_j)<\epsilon\\
0 & \textrm{if} & d(f(p_i),p_j)\geq\epsilon\end{array}\right.$$
and a surjective map $\beta:\Sigma_A\rightarrow \tilde{\Lambda}$.

\begin{cla}\label{c.periodic} If
$\{p_i\}_{i=1}^{N}$ consists of the orbits of periodic points,
then $\overline{\mathrm{Per}(\tilde{\Lambda})}=\tilde{\Lambda}$.
\end{cla}

\noindent{\bf Proof.}  It is sufficient to show periodic points
are dense in $\Sigma_A$.  Let $s\in\Sigma_A$ and $s[-n,n]$ be a
word of length $2n$ in $s$.  Fix $i\in[-n,n-1]$.  Then
$a_{s_is_{i+1}}=1$.  To complete the proof we show there is a word
$w_i$ in $\Sigma_A$ starting with $s_{i+1}$ and ending with
$s_{i}$.

Let $m_{i+1}$ be the period of $p_{s_{i+1}}$ and $m_i$ be the
period of $p_{m_i}$.  Take the sequence of points
$$p_{s_{i+1}},
f(p_{s_{i+1}}),...,f^{m_{i+1}-1}(p_{s_{i+1}}),f(p_{s_i}),...,f^{m_i-1}(p_{s_i}),f^{m}(p_{s_i})=p_{s_i},$$
and $w_i$ to be the associated word in $\Sigma_A$.

Then there is a periodic point $t\in\Sigma_A$ containing the word
$$s[-n,n]w_{n-1}\cdots w_{-n}.$$  Hence, periodic points are dense
in $\Sigma_A$. $\Box$

\begin{lem}\label{l.periodic}
Let $\Lambda$ be a non-locally maximal hyperbolic set with periodic
points dense and $V$ be a neighborhood of $\Lambda$.  Then there
exists a hyperbolic set $\tilde{\Lambda}\supset\Lambda$ with a
Markov partition contained in the closure of $V$ and a sequence of
periodic points $p_n\in\tilde{\Lambda}-\Lambda$ converging to a
point $x\in\Lambda$.
\end{lem}

\noindent{\bf Proof.}  Let $\nu>0$ and $\eta>0$ be defined as in
the proof of Theorem 1.5 in~\cite{Fis1} so that associated to any $\nu$ dense set of points in $\Lambda$ there is a subshift of finite type associated with the transition matrix as described above, and so that for any two points $x,y\in\Lambda$, where $d(x,y)<\nu$, the intersection $W^s_{\eta}(x)\cap W^u_{\eta}(y)$ consists of one point in $\overline{V}$.  

If $\Lambda$ is a
non-locally maximal hyperbolic set, then $\Lambda$ does not have a
local product structure.  Hence, there exists points
$x,y\in\Lambda$ such that $d(x,y)<\nu/2$ and $W^s_{\eta}(x)\cap
W^u_{\eta}(y)$ is not contained in $\Lambda$.  Since periodic
points are dense in $\Lambda$ the continuity of stable and
unstable manifolds implies there exist periodic points $q_1$ and
$q_2$ in $\Lambda$, of period $m_1$ and $m_2$, respectively, such
that $d(q_1,q_2)<\nu$ and $W^s_{\eta}(q_1)\cap W^u_{\eta}(q_2)$ is
not contained in $\Lambda$.

Let $\{p_i\}_{i=1}^N$ be a dense set of $\nu$ dense periodic
points containing the orbit of each periodic point, and containing
the points $q_1$ and $q_2$.  Then there are words $w_1$ and $w_2$
in $\Sigma_A$ corresponding to the orbits
$q_1,f(q_1),...,f^{m_1-1}(q_1)$ and
$q_2,f(q_2),...,f^{m_2-1}(q_2)$, respectively.  Let
$s=(\overline{w_2}.\overline{w_1})\in\Sigma_A$.  Then
$$z=\beta(s)=W^s_{\eta}(q_1)\cap W^u_{\eta}(q_2)\in\tilde{\Lambda}.$$

Since $\tilde{\Lambda}$ is an invariant set we know that
$\mathcal{O}(z)\subset\tilde{\Lambda}-\Lambda$.  Furthermore, since $z\in W^s_{\eta}(q_1)$ and $q_1\in\Lambda$ we know that  $f^n(z)$
converges to $\Lambda$ as $n\rightarrow\infty$.  

From
Claim~\ref{c.periodic} the periodic points of $\tilde{\Lambda}$
are dense.  Hence there is a sequence of periodic points
$p_n\in\tilde{\Lambda}-\Lambda$ such that $d(p_n,f^{nm_1}(z))<\nu/2^n$.  Therefore, $p_n$ converges to $q_1$.    $\Box$
\\

\noindent{\bf Proof of Theorem~\ref{t.hypcommute}} We first note
that from Proposition~\ref{p.abcd} there exists a neighborhood $U$ of $\Lambda$ such that if $p$
is a periodic point of $f$ and $p\in U$, then $p\in\Lambda$.

Suppose $\Lambda$ is not locally maximal for $g$.  Then there
exists a neighborhood $V$ of $\Lambda$ contained in $U$ such that
$$\Lambda_V=\bigcap_{n\in\mathbb{Z}}f^n(\overline{V})\supsetneq\Lambda$$
is a hyperbolic set.  From Lemma~\ref{l.periodic} there exists a
hyperbolic set $\tilde{\Lambda}\supset\Lambda$ with a Markov
partition that is contained in $\Lambda_V$ and a sequence of
periodic points $p_n\in\tilde{\Lambda}-\Lambda$ converging to a
point $x\in\Lambda$.  Since each $p_n\in\mathrm{Per}_h(g)$ we know the points $p_n$
are periodic points for $f$ contained in $U-\Lambda$, a
contradiction. Hence, $\Lambda$ is locally maximal for $g$. $\Box$

\bibliographystyle{plain}

\end{document}